\documentclass[11pt]{amsart}

\newcommand{\qbin}[2]{\genfrac{[}{]}{0pt}{}{#1}{#2}}
\newcommand{\qbins}[2]{{\textstyle\genfrac{[}{]}{0pt}{}{#1}{#2}}}

\begin{document}

\title[$q$-Hypergeometric proofs]{$q$-Hypergeometric proofs of
polynomial analogues of the triple product identity, Lebesgue's identity 
and Euler's pentagonal number theorem}

\author[Ole Warnaar]{S. Ole Warnaar}

\address{Department of Mathematics and Statistics,
The University of Melbourne, Vic 3010, Australia}
\email{warnaar@ms.unimelb.edu.au}


\thanks{Work supported by the Australian Research Council}

\begin{abstract}
We present alternative, $q$-hypergeometric proofs of some
polynomial analogues of classical $q$-series identities
recently discovered by Alladi and Berkovich, and Berkovich and Garvan.
\end{abstract}

\maketitle

\section{Introduction}
In two recent papers, Alladi and Berkovich \cite{AB02} and
Berkovich and Garvan \cite{BG02} proved the following three
polynomial identities:
\begin{multline}\label{id1}
\sum_{n=0}^L \frac{z^{-n}+z^{1+n}}{1+z}\, q^{T_n} \\ =
\sum_{i,j,k\geq 0}(-1)^k z^{i-j}q^{T_i+T_j+T_k}
\qbin{L-i}{j}\qbin{L-j}{k}\qbin{L-k}{i},
\end{multline}
\begin{multline}\label{id2}
\sum_{i,j\ge 0}(-1)^j z^{2j}q^{T_i+T_j}\qbin{L-j}{i}\qbin{i}{j} \\
=\sum_{i,j,k\geq 0} (-1)^j z^{i+j}q^{T_i+T_j+T_k}
\qbin{L-i}{j}\qbin{L-j}{k}\qbin{L-k}{i}
\end{multline}
and
\begin{equation}\label{id3}
\sum_{j=-\infty}^{\infty} (-1)^j q^{j(3j+1)/2}\qbin{2L-j}{L+j}=1.
\end{equation}
Here $T_n=n(n+1)/2$ is a triangular number and $\qbins{n}{k}$ a 
$q$-binomial coefficient defined as
\begin{equation*}
\qbin{n}{k}_q=\qbin{n}{k}=\frac{(q;q)_n}{(q;q)_k(q;q)_{n-k}}
\end{equation*}
for $0\leq k\leq n$ and zero otherwise, with
$(a;q)_n=\prod_{k=0}^{n-1}(1-aq^k)$.

Identity \eqref{id1}, which is \cite[Eq. (1.12)]{AB02} and which was
first stated in \cite{AB01}, is a polynomial analogue of Jacobi's
triple product identity.
Identity \eqref{id2}, which is \cite[Eq. (1.15)]{AB02}, is a polynomial
analogue of Lebesgue's identity and \eqref{id3}, which is 
\cite[Eq. (1.28)]{BG02}, is a polynomial analogue of Euler's pentagonal 
number theorem.

The first explicit proof of \eqref{id1} was given by Berkovich and Riese 
\cite{BR02} who showed that both sides satisfy a highly nontrivial
fourth order recurrence relation. The proof of \eqref{id1} given by
Alladi and Berkovich uses transformation formulas for basic hypergeometric
series. In particular they show that the right side of \eqref{id1} can be
transformed into the left side using, consecutively, Heine's ${_2\phi_1}$
transformation, the $q$-Chu--Vandermonde sum and the Sears--Carlitz 
transformation between ${_3\phi_2}$ and ${_5\phi_4}$ series.
In this note we present a simple, one-page proof of \eqref{id1} 
that only requires elementary summations and no transformations.
As a bonus we find that \eqref{id1} is the $c=0$ instance of
\begin{multline}\label{id1b}
\sum_{n=-L}^L \sum_{m=-L}^n
(-1)^{n+m} z^m q^{T_n}\frac{(cq;q)_{L-m}}{(cq;q)_{L-n}(cq;q)_{n-m}} \\ =
\sum_{i,j,k\geq 0} (-1)^k z^{i-j} q^{T_i+T_j+T_k}\frac{1-c}{1-cq^j}
\qbin{L-i}{j}\qbin{L-j}{k}\qbin{L-k}{i}
\end{multline}
and of
\begin{multline}\label{id1c}
\sum_{n=0}^L \sum_{m=-n}^n (-1)^{n+m} z^m q^{T_n}
\frac{(cq;q)_{L-|m|}}{(cq;q)_{L-n}(cq;q)_{n-|m|}} \\=
\sum_{i,j,k\geq 0} (-1)^k z^{i-j} q^{T_i+T_j+T_k}
\frac{1-c}{1-cq^{\min(i,j)}}
\qbin{L-i}{j}\qbin{L-j}{k}\qbin{L-k}{i}.
\end{multline}
Letting $c$ tend to infinity and performing the sum over $m$ leads to the
following two variants of \eqref{id1}:
\begin{multline*}
\sum_{n=0}^{2L}
\frac{(-1)^n z^{-L}q^{n(2L-n+1)} +z^{L-n+1}}{q^n+z}\, q^{T_{L-n}} \\
= \sum_{i,j,k\geq 0} (-1)^k z^{i-j} q^{T_i+T_j+T_k-j}
\qbin{L-i}{j}\qbin{L-j}{k}\qbin{L-k}{i}
\end{multline*}
and
\begin{multline*}
\sum_{n=0}^L \biggl\{
\frac{(-1)^n q^{(n+1)(L-n)}+z^{n+1}}{q^{L-n}+z}
+\frac{(-1)^{n+1} q^{n(L-n)}+z^{-n}}{1+zq^{L-n}}
\biggr\}\, q^{T_n} \\
=\sum_{i,j,k\geq 0} (-1)^k z^{i-j} q^{T_i+T_j+T_k-\min(i,j)}
\qbin{L-k}{i}\qbin{L-j}{k}\qbin{L-i}{j}.
\end{multline*}

The proof of \eqref{id2} as found by Alladi and Berkovich \cite{AB02} uses 
combinatorial methods. They however ask for a $q$-hypergeometric proof
of \eqref{id2}, or more precisely, for a $q$-hypergeometric proof of the
closely related \cite[Eq. (6.14)]{AB02}
\begin{multline}\label{id2b}
\sum_{j\ge 0}(z^2 q^2;q^2)_j q^{\binom{L-2j}{2}}\qbin{L+1}{2j+1} \\
=\sum_{i,j,k\geq 0} (-1)^j z^{i+j}q^{T_i+T_j+T_k}
\qbin{L-i}{j}\qbin{L-j}{k}\qbin{L-k}{i}.
\end{multline}
Here we present such a $q$-hypergeometric proof of \eqref{id2} and 
\eqref{id2b}.

Identity \eqref{id3} was discovered and proved by Berkovich and Garvan by
finitizing Dyson's proof of Euler's identity using the rank of a partition.
Again it is posed as a problem to find a $q$-hypergeometric proof.
Here we show that \eqref{id3} is a simple consequence of a cubic 
summation formula of \cite{W00}. This same cubic sum may be applied to 
also yield
\begin{equation}\label{id4}
\sum_{j=-\infty}^{\infty} (-1)^j q^{j(3j-1)/2}\qbin{2L-j+1}{L+j}=1.
\end{equation}

\section{Proofs}
In the following we adopt the notation of \cite{GR90} for basic 
hypergeometric series, writing
\begin{multline*}
{_{r+1}\phi_r}\biggl[\genfrac{}{}{0pt}{}{a_1,a_2,\dots,a_{r+1}}
{b_1,\dots,b_r};q,z\biggr]=
{_{r+1}\phi_r}(a_1,a_2,\dots,a_{r+1};b_1,\dots,b_r;q,z) \\
=\sum_{k=0}^{\infty}\frac{(a_1,a_2,\dots,a_{r+1};q)_k}
{(q,b_1,\dots,b_r;q)_k}\,z^k,
\end{multline*}
where $(a_1,a_2,\dots,a_k;q)_n=(a_1;q)_n\cdots(a_k;q)_n$.
Moreover, whenever a basic hypergeometric identity occurs
with a term $q^{-n}$ and/or $q^{-m}$, $n$ and $m$ are assumed to be
nonnegative integers.

\subsection{Proof of \eqref{id1}}
Using $(1+z^{2n+1})/(1+z)=\sum_{m=0}^{2n}(-z)^m$, interchanging sums
and shifting $n\to n+|m|$, the left side of \eqref{id1} can be written as
\begin{equation}\label{zexp}
\text{LHS}\eqref{id1}=
\sum_{m=-L}^L \sum_{n=0}^{L-|m|} (-1)^{n} z^m q^{T_{n+|m|}}.
\end{equation}
If on the right of \eqref{id1} we eliminate $i$ in favour of
$m$ by $i=m+j$ we can use \eqref{zexp} to equate coefficients of 
$z^m$ in \eqref{id1}.
This leads to an identity, denoted $(\ast)$, which should hold
for all $m$ and $L$ such that $|m|\leq L$.
After dividing both sides by $q^{T_{|m|}}$, $(\ast)$ can easily be
seen to be the $b=q^{|m|+1}$ and $N=L-|m|$ instance of
\begin{equation}\label{bN}
\sum_{n=0}^N (-1)^n b^n q^{\binom{n}{2}} =
\sum_{j=0}^{\lfloor N/2\rfloor}\sum_{k=0}^{N-j}
\frac{(-1)^k b^j q^{j^2+T_k}(q,b;q)_{N-j}(b;q)_{N-k}}
{(q,b;q)_j(q,b;q)_{N-j-k}(q;q)_k(q;q)_{N-2j}}.
\end{equation}
We now take the right-hand side and make the variable change
$k\to N-j-k$. Then
\begin{align*}
\text{RHS}\eqref{bN}&=
\sum_{j=0}^{\lfloor N/2\rfloor}
\frac{(-1)^{N-j} b^j q^{j^2+T_{N-j}}(b;q)_{N-j}}{(q;q)_j(q;q)_{N-2j}}\,
{_2\phi_1}\biggl[\genfrac{}{}{0pt}{}{q^{-N+j},bq^j}{b};q,1\biggr] \\
&=\sum_{j=0}^{\lfloor N/2\rfloor}
b^j q^{j^2}(bq^j;q)_{N-2j}\qbin{N-j}{j} \\
&=\sum_{j=0}^{\lfloor N/2\rfloor}\sum_{n=0}^{N-2j}
(-1)^n b^{j+n}q^{j(j+n)+\binom{n}{2}}\qbin{N-j}{j}\qbin{N-2j}{n} \\
&=\sum_{n=0}^N (-1)^n b^n q^{\binom{n}{2}}\qbin{N}{n}
{_2\phi_1}\biggl[\genfrac{}{}{0pt}{}{q^{-n},q^{-N+n}}{q^{-N}};q,q\biggr] \\
&=\sum_{n=0}^N (-1)^n b^n q^{\binom{n}{2}}\qbin{N}{n}
\frac{(q^{-n};q)_n}{(q^{-N};q)_n}q^{-(N-n)n}=\text{LHS}\eqref{bN}.
\end{align*}
Here the second equality follows from
\begin{equation}\label{z1}
{_2}\phi_1(a,q^{-n};aq^{-m};q,1)=\frac{(q^{-n};q)_n}{(aq^{-m};q)_m},
\end{equation}
the third equality follows from the $q$-binomial theorem 
\cite[Eq. (II.4)]{GR90}
\begin{equation}\label{qbthm}
{_1\phi_0}(q^{-n};\text{---}\,;q,z)=(zq^{-n};q)_n,
\end{equation}
the fourth equality follows from the shift $n\to n-j$ followed by 
an interchange of sums, and the second-last equality follows from 
the $q$-Chu--Vandermonde sum \cite[Eq. (II.6)]{GR90}
\begin{equation}\label{qCV}
{_2\phi_1}(a,q^{-n};c;q,q)=\frac{(c/a;q)_n}{(c;q)_n}\, a^n.
\end{equation}
To prove \eqref{z1}, take Heine's transformation \cite[Eq. (III.1)]{GR90}
\begin{equation*}
{_2}\phi_1(a,b;c;q,z)=\frac{(a,bz;q)_{\infty}}{(c,z;q)_{\infty}}\,
{_2}\phi_1(c/a,z;bz;q,a),
\end{equation*}
specialize $b=q^{-n}$ and $c=aq^{-m}$ and use 
$(aq^{-n};q)_{\infty}/(a;q)_{\infty}=(aq^{-n};q)_n$. After letting $z$ 
tend to $1$ this yields \eqref{z1}.

By replacing the $q$-Chu--Vandermonde sum \eqref{qCV} with the more general 
$q$-Saalsch\"utz sum \cite[Eq. (II.12)]{GR90}, the above proof generalizes 
to yield 
\begin{multline*}
\sum_{n=0}^N (-1)^n b^n q^{\binom{n}{2}}
\frac{(cq;q)_N}{(cq;q)_n(cq;q)_{N-n}} \\=
\sum_{j=0}^{\lfloor N/2\rfloor}\sum_{k=0}^{N-j}
\frac{(-1)^k b^j q^{j^2+T_k} 
(c;q)_j(q,b;q)_{N-j}(b;q)_{N-k}}
{(q,b,cq;q)_j(q,b;q)_{N-j-k}(q;q)_k(q;q)_{N-2j}},
\end{multline*}
which for $c=0$ reduces to \eqref{bN}. Taking $b=q^{m+1}$ 
($b=q^{|m|+1}$) and $N=L-m$ ($N=L-|m|$), multiplying
both sides by $b^m$ and summing $m$ from $-L$ to $L$ gives \eqref{id1b}
(\eqref{id1c}).

\subsection{Proof of \eqref{id2} and \eqref{id2b}}
First consider the left-hand side of \eqref{id2}. Shifting $i\to i+j$ 
this becomes
\begin{align}\label{LHS}
\text{LHS}\eqref{id2}&=\sum_{j=0}^{\lfloor L/2 \rfloor}
(-1)^j z^{2j} q^{j(j+1)}\qbin{L-j}{j}
{_{1}\phi_0}(q^{-L+2j};\text{---}\,;q,-q^{L-j+1}) \\
&=\sum_{j=0}^{\lfloor L/2 \rfloor}
(-1)^j z^{2j} q^{j(j+1)}\qbin{L-j}{j}\frac{(-q;q)_{L-j}}{(-q;q)_j}, \notag
\end{align}
where the second equality follows from \eqref{qbthm}.

Next consider the right-hand side of \eqref{id2}.
By reshuffling the terms that make up the three $q$-binomial coefficients,
it readily follows that the summand $S_{L;i,j,k}$ on the right
satisfies $S_{L;i,j,k}=(-1)^{i+j}S_{L;j,i,k}$. Hence all contributions to 
the sum arising from $i+j$ being odd cancel, and we may add the 
restriction ``$i+j$ even'' to the sum.
By replacing $j\to 2j-i$ we can then extract the coefficient of $z^{2j}$.
Equating this with the coefficient of $z^{2j}$ arising from \eqref{LHS}
results in the identity
\begin{equation}\label{z2j}
\sum_{i,k\geq 0} \frac{(-1)^{i-j} q^{(i-j)^2+T_k}}{(-q;q)_{L-2j}}
\qbin{L-i}{2j-i}\qbin{L+i-2j}{k}\qbin{L-k}{i}
=\qbin{L-j}{j}_{q^2}
\end{equation}
for $0\leq j\leq \lfloor L/2\rfloor$. This is easily proved as follows:
\begin{align*}
\text{LHS}\eqref{z2j}&=
\qbin{L}{2j}\sum_{i=0}^{2j}\frac{(-1)^{i-j} q^{(i-j)^2}}
{(-q;q)_{L-2j}}\qbin{2j}{i}
{_2\phi_1}\biggl[\genfrac{}{}{0pt}{}{q^{-L+i},q^{-L-i+2j}}
{q^{-L}};q,-q^{L-2j+1}\biggr]\\
&=\qbin{L}{2j}\sum_{i=0}^{2j}(-1)^{i-j} q^{(i-j)^2} \qbin{2j}{i}
{_2\phi_1}\biggl[\genfrac{}{}{0pt}{}{q^{-i},q^{-2j+i}}
{q^{-L}};q,-q\biggr]\\
&=(-1)^j q^{j^2}\qbin{L}{2j}
\sum_{k=0}^j\frac{(q^{-2j};q)_{2k}\, q^k}{(q,q^{-L};q)_k}
\sum_{i=0}^{2j-2k}\frac{q^{T_i} (q^{-2j+2k};q)_i}{(q;q)_i} \\
&=\qbin{L}{2j}
\sum_{k=0}^j (-1)^k q^{k(2j-k+1)}\frac{(q^{-2j};q)_{2k}(q;q^2)_{j-k}}
{(q,q^{-L};q)_k} \\[1.5mm]
&=(q;q^2)_j\qbin{L}{2j}{_2\phi_1}\biggl[\genfrac{}{}{0pt}{}{q^{-j},-q^{-j}}
{q^{-L}};q,q\biggr]\\[2mm]
&=(q;q^2)_j\qbin{L}{2j}
\frac{(-q^{-L+j};q)_j}{(q^{-L};q)_j} (-1)^j q^{-j^2}=\text{RHS}\eqref{z2j}.
\end{align*}
Here the second equality follows from Heine's transformation 
\cite[Eq. (III.3)]{GR90}
\begin{equation*}
{_2}\phi_1(a,b;c;q,z)=\frac{(abz/c;q)_{\infty}}{(z;q)_{\infty}}\,
{_2}\phi_1(c/a,c/b;c;q,abz/c),
\end{equation*}
the third equality follows from an interchange sums and the shift 
$i\to i+k$, the fourth equality follows from 
\begin{equation}\label{ex26}
\sum_{i=0}^{2n}\frac{q^{T_i} (q^{-2n};q)_i}{(q;q)_i}
=(-1)^n q^{-n^2}(q;q^2)_n
\end{equation}
and the second-last equality follows from \eqref{qCV}.
To prove \eqref{ex26}, take \cite[Excer. 2.6]{GR90}
\begin{equation*}
{_4\phi_3}\Bigl[\genfrac{}{}{0pt}{}{q^{-n},b,c,-q^{1-n}/bc}
{q^{1-n}/b,q^{1-n}/c,-bc};q,q\Bigr]
=\begin{cases}\displaystyle
\frac{(q,b^2,c^2,q^2)_{n/2}(bc;q)_n}{(b,c,q)_n(b^2c^2;q^2)_{n/2}}
& n\text{ even} \\[3mm]
0 & n\text{ odd},
\end{cases}
\end{equation*}
replace $n\to 2n$ and let $b$ and $c$ tend to infinity.

To now also prove \eqref{id2b} we only need to show its left-hand side
equals the left-hand side of \eqref{id2}.
To achieve this we expand $(z^2 q^2;q^2)_j$ using \eqref{qbthm},
expressing the left side of \eqref{id2b} as a double sum over $j$ and $k$.
Interchanging the order of these sums and shifting $j\to j+k$ leads to
\begin{multline*}
\text{LHS}\eqref{id2b}=\sum_{k=0}^{\lfloor L/2\rfloor}
(-1)^k z^{2k} q^{k(k+1)+\binom{L-2k}{2}}\qbin{L+1}{2k+1} \\
\times {_2\phi_1}\biggl[\genfrac{}{}{0pt}{}{q^{-L+2k},q^{1-L+2k}}
{q^{2k+3}};q^2,q^2\biggr].
\end{multline*}
The ${_2\phi_1}$ can be evaluated by \eqref{qCV} with $n=L/2-k$ 
when $L$ is even and $n=(L-1)/2-k$ when $L$ is odd.
Specifically, if we set $L=2N+\sigma$ with 
$\sigma\in\{0,1\}$ then the above ${_2\phi_1}$ yields
\begin{equation*}
\frac{(q^{2N+2\sigma+2};q^2)_{N-k}}{(q^{2k+3};q^2)_{N-k}}\,
q^{-\binom{L-2k}{2}}=\frac{(q^2;q^2)_{L-k}(q;q)_{2k+1}}
{(q;q)_{L+1}(q^2;q^2)_k}\, q^{-\binom{L-2k}{2}},
\end{equation*}
which is independent of the parity of $L$. Hence
\begin{equation*}
\text{LHS}\eqref{id2b}=\sum_{k=0}^{\lfloor L/2\rfloor}
(-1)^k z^{2k} q^{k(k+1)}\frac{(q^2;q^2)_{L-k}}{(q^2;q^2)_k(q;q)_{L-2k}}.
\end{equation*}
Since this is identical to the expression for the left-hand side 
of \eqref{id2} found in \eqref{LHS} we are done.

\subsection{Proof of \eqref{id3} and \eqref{id4}}
In \cite[Cor. 4.13]{W00} we proved (an elliptic analogue of)
the cubic summation
\begin{multline*}
\sum_{k=0}^{\lfloor n/2 \rfloor} \frac{1-aq^{4k}}{1-a}
\frac{(a,aq^{n+1};q^3)_k}{(q,q^{-n};q)_k}
\frac{(q^{-n};q)_{2k}}{(aq^{n+1};q)_{2k}}
\frac{(c,d;q)_k}{(aq^3/c,aq^3/d;q^3)_k}\, q^k \\
=\begin{cases}\displaystyle
\frac{(aq^3,q^{2-n}/c,q^{2-n}/d;q^3)_{\lfloor n/3 \rfloor}}
{(aq^3/c,aq^3/d,q^{2-n}/cd;q^3)_{\lfloor n/3 \rfloor}}
& n\not\equiv 2\pmod{3} \\[3mm]
0 & n\equiv 2\pmod{3}
\end{cases}
\end{multline*}
for $cd=aq^{n+1}$.
Replacing $a\to a^2$ followed by $c\to ac$ and $d\to ad$, and then 
letting $a$ tend to zero yields
\begin{equation*}
\sum_{k=0}^{\lfloor n/2 \rfloor} 
\frac{(q^{-n};q)_{2k}}{(q,q^{-n};q)_k}\, q^k
=\begin{cases}(-1)^{\lfloor n/3 \rfloor}q^{-n(n-1)/6}
& n\not\equiv 2\pmod{3} \\[3mm]
0 & n\equiv 2\pmod{3}.
\end{cases}
\end{equation*}
Taking $n=3L$, replacing $k\to j+L$ and making some simplifications
gives \eqref{id3}. In much the same way does $n=3L+1$ lead to \eqref{id4}.
Making the same replacements for $a,b$ and $c$ but letting $a$ tend to 
infinity instead of zero results in
\begin{equation*}
\sum_{k=0}^{\lfloor n/2 \rfloor}(-1)^k q^{\binom{k}{2}} \qbin{n-k}{k}
=\begin{cases}(-1)^{\lfloor n/3 \rfloor}q^{n(n-1)/6}
& n\not\equiv 2\pmod{3} \\[3mm]
0 & n\equiv 2\pmod{3}.
\end{cases}
\end{equation*}

\end{document}